\documentclass[12pt]{article}
\def\bC{\mathbf{\overline{C}}}
\begin{document}
\def\supp{\mathrm{supp}\,}
\def\C{\mathbf{C}}
\def\Z{\mathbf{Z}}
\def\Del{\nabla}
\def\P{\mathbf{P}}
\def\Rea{{\mathrm{Re}\,}}
\author{Matthew Barrett and Alexandre Eremenko\thanks{Supported by
NSF grant DMS-055279 and by the Humboldt Foundation.}}
\title{Generalization of a theorem of Clunie and Hayman}
\maketitle
\begin{abstract} Clunie and Hayman proved that
if the spherical derivative $\| f'\|$ of
an entire function satisfies $\| f'\|(z)=O(|z|^\sigma)$
then $T(r,f)=O(r^{\sigma+1}).$
We generalize this to holomorphic curves
in projective space of dimension $n$
omitting $n$ hyperplanes in general position.

MSC 32Q99, 30D15.
\end{abstract}

\noindent
{\bf Introduction}
\vspace{.1in}

We consider holomorphic curves $f:\C\to\P^n$; for the
general background on the subject we refer to \cite{Lang}.
The Fubini--Study derivative $\| f'\|$ measures the
length distortion from the
Euclidean metric in $\C$ to the
Fubini--Study metric in $\P^n$. The explicit expression
is 
$$\| f'\|^2=
\| f\|^{-4}\sum_{i<j}|f_i^\prime f_j-f_if_j^\prime|^2,$$
where $(f_0,\ldots,f_n)$ is a homogeneous
representation of $f$
(that is the $f_j$ are entire functions which
never simultaneously
vanish), and
$$\| f\|^2=\sum_{j=0}^n|f_j|^2.$$
See \cite{CE} for a general discussion of the
Fubini-Study derivative.

We recall that the Nevanlinna--Cartan characteristic is
defined by
$$T(r,f)=\int_0^r \frac{dt}{t}
\left(\frac{1}{\pi}\int_{|z|\leq t}\| f'\|^2(z)dm(z)\right),$$
where $dm$ is the area element in $\C$.
So the condition
\begin{equation}
\label{cond}
\limsup_{z\to\infty} |z|^{-\sigma}\| f'(z)\|\leq K<\infty
\end{equation}
implies
\begin{equation}
\label{1}
\limsup_{r\to\infty}\frac{T(r,f)}{r^{2\sigma+2}}<\infty.
\end{equation}

Clunie and Hayman \cite{CH} found that for 
curves $\C\to\P^1$ omitting one
point in $\P^1$, a stronger conclusion follows from
(\ref{cond}), namely 
\begin{equation}
\label{concl}
\limsup_{r\to\infty}\frac{T(r,f)}{r^{\sigma+1}}
\leq KC(\sigma).
\end{equation}
In the most important case $\sigma=0$, a
different proof of this fact for $n=1$ is due to 
Pommerenke \cite{Pom}. Pommerenke's method gives the
exact constant $C(0)$.
In this paper we prove that this phenomenon
persists in all dimensions.
\vspace{.1in}

\noindent
{\bf Theorem.} {\em For holomorphic curves
$f:\C\to\P^n$
omitting $n$ hyperplanes in general position, condition
$(\ref{cond})$ implies $(\ref{concl})$ with an explicit
constant $C(n,\sigma)$.}
\vspace{.1in}

In \cite{E}, the case $\sigma=0$ was considered.
There it was proved that holomorphic curves in $\P^n$
with
bounded spherical derivative and omitting $n$
hyperplanes in general position must satisfy
$T(r,f)=O(r)$.
With a stronger assumption that $f$ omits $n+1$ hyperplanes
this was earlier established by
Berteloot and Duval \cite{BD} and 
by Tsukamoto \cite{T2}.
The proof in \cite{E} has two drawbacks: it does not
extend to arbitrary $\sigma\geq 0$,
and it is non-constructive; unlike
Clunie--Hayman and Pommerenke's proofs mentioned above,
it does not give an explicit constant in (\ref{concl}).

It is shown in \cite{E} that the condition
that $n$ hyperplanes are omitted is exact:
there are curves in any dimension $n$ satisfying
(\ref{cond}), $T(r,f)\sim cr^{2\sigma+2}$ and omitting $n-1$ hyperplanes.
\vspace{.2in}

\noindent
{\bf Preliminaries}\nopagebreak
\vspace{.2in}

Without loss of generality
we assume that
the omitted hyperplanes are given in the
homogeneous coordinates
by the equations
$\{ w_j=0\},\; 1\leq j\leq n.$
We fix a homogeneous representation
$(f_0,\ldots,f_n)$ of our curve, where $f_j$
are entire functions, and $f_n=1$.
Then
\begin{equation}
\label{u}
u=\log\sqrt{|f_0|^2+\ldots+|f_n|^2}
\end{equation}
is a positive subharmonic function, and Jensen's formula
gives
$$T(r,f)=\frac{1}{2\pi}\int_{-\pi}^\pi u(re^{i\theta})d\theta-u(0)=\int_0^r\frac{n(t)}{t}dt,$$
where $n(t)=\mu(\{ z:|z|\leq t\})$, and $\mu=\mu_u$
is the Riesz measure  of $u$, that is the measure with
the density
\begin{equation}\label{3}
\frac{1}{2\pi}\Delta u=\frac{1}{\pi}\| f'\|^2.
\end{equation}
This measure $\mu$ is also called Cartan's measure of $f$.
Positivity of $u$ and (\ref{1}) imply that
all $f_j$ are of order at most $2\sigma+2$, normal type.
As $f_j(z)\neq 0,\; 1\leq j\leq n$ we conclude that 
$$f_j=e^{P_j},\quad 1\leq j\leq n,$$
where
\begin{equation}\label{starr}
P_j\quad\mbox{are polynomials of degree at most}
\quad 2\sigma+2.
\end{equation}

We need two lemmas from potential theory.
\vspace{.1in}

\noindent
{\bf Lemma 1.}$\;$\cite{E} {\em Let $v$ be a non-negative
harmonic function in the closure of the
disc $B(a,R)$, and assume that
$v(z_1)=0$ for some point $z_1\in\partial B(a,R).$
Then} 
$$v(a)\leq2R|\Del v(z_1)|.$$
\vspace{.1in}

\noindent
{\bf Lemma 2.} {\em Let $v$ be a non-negative superharmonic function
in the closure of the disc $B(a,R)$, and suppose that
$v(z_1)=0$ for some $z_1\in\partial B(a,R)$.
Then
$$\left|\mu_v(B(a,R/2))\right|\leq 3R|\Del v(z_1)|.$$}
\vspace{.1in}

{\em Proof.} Function $v(a+Rz)$ satisfies the conditions
of the lemma with $R=1$. So it is enough to prove the lemma with
$a=0$ and $R=1$. Let
$$w(z)=\int_{|\zeta|<1}G(z,\zeta)d\mu_v(\zeta)$$
be the Green potential of $\mu_v$. Then $w\leq v$ and $w(z_1)=v(z_1)$
which implies that 
$$|\Del v(z_1)|\geq\left|\frac{\partial w}{\partial|z|}(z_1)\right|.$$
Minimizing $|\partial G/\partial|z||$ over $|z|=1$ and $|\zeta|=1/2$
we obtain $1/3$ which proves the lemma.
\vspace{.2in}

\noindent
{\bf Proof of the theorem}
\vspace{.2in}

We may assume without loss of generality
that $f_0$ has infinitely many zeros.
Indeed, we can compose $f$ with an automorphism
of $\P^n$, for example replace $f_0$ by $f_0+cf_1,\; c\in\C$
and leave all other $f_j$ unchanged.
This transformation changes neither the $n$
omitted hyperplanes nor the rate of growth of $T(r,f)$
and multiplies the spherical derivative by
a bounded factor. 

Put $u_j=\log|f_j|$, and 
$$u^*=\max_{1\leq j\leq n}u_j.$$
Here and in what follows $\max$ denotes the
pointwise maximum of subharmonic functions. 
\vspace{.1in}

\noindent
{\bf Proposition 1.} {\em Suppose that at some point $z_1$
we have $$u_m(z_1)=u_k(z_1)\geq u_j(z_1)$$
for some $m\neq k$ and all $j$; $m,k,j\in\{0,\ldots,n\}$.
Then}
$$\| f'(z_1)\|\geq (n+1)^{-1}
|\nabla u_m(z_1)-\nabla u_k(z_1)|.$$

{\em Proof.}
$$\| f'(z_1)\|\geq
\frac{|f_m^\prime(z_1)f_k(z_1)
-f_m(z_1)f_k^\prime(z_1)|}{|f_0(z_1)|^2+
\ldots+|f_n(z_1)|^2}
\geq(n+1)^{-1}\left|\frac{f_m^\prime(z_1)}{f_m(z_1)}-
\frac{f_k^\prime(z_1)}{f_k(z_1)}\right|,$$
and the conclusion of the proposition follows
since
$|\nabla\log|f||=|f'/f|$.
\vspace{.1in}

\noindent
{\bf Proposition 2.} {\em For every $\epsilon>0$, we have
$$u(z)\leq u^*(z)+K(2+\epsilon)^{\sigma+1} (n+1)|z|^{\sigma+1}$$
for all
$|z|>r_0(\epsilon)$.}
\vspace{.1in}

{\em Proof.} If $u_0(z)\leq u^*(z)$ for all sufficiently
large $|z|$, then there is nothing to prove.
Suppose that $u_0(a)>u^*(a)$, and consider the largest disc $B(a,R)$
centered at $a$ where the inequality $u_0(z)>u^*(z)$ persists.
If $z_0$ is the zero of the
smallest modulus of $f_0$ then
$R\leq |a|+|z_0|<(1+\epsilon)|a|$ when $|a|$ is large enough.

There is a point $z_1\in\partial B(a,R)$ such that $u_0(z_1)=u^*(z_1)$.
This means that there is some $k\in\{1,\ldots,n\}$ such that
$u_0(z_1)=u_k(z_1)\geq u_m(z_1)$ for all $m\in\{1,\ldots,n\}$.
Applying Proposition 1 we obtain
$$|\nabla u_k(z_1)-\nabla u_0(z_1)|\leq (n+1) \| f'(z_1)\|.$$
Now $u_0(z)>u^*(z)\geq u_k(z)$ for $z\in B(a,R)$, so we can apply Lemma 1
to $v=u_0-u_k$ in the disc $B(a,R)$.
This gives
$$u_0(a)-u_k(a)\leq 2R|\nabla u_k(z_1)-\nabla u_0(z_1)|
\leq 2R(n+1) \| f'(z_1)\|.$$
Now $R<(1+\epsilon)|a|$ and $|z_1|\leq (2+\epsilon)|a|$, so
$$u_0(a)\leq u^*(a)+
K(2+\epsilon)^{\sigma+1} (n+1)|a|^{\sigma+1}, $$
and the result follows because $u=\max\{ u_0,u^*\}+O(1)$.
\vspace{.1in}

Next we study the Riesz measure of the subharmonic
function $$u^*=\max\{u_1,\ldots,u_n\}.$$
We begin with maximum of two harmonic functions.
Let $u_1$ and $u_2$ be two harmonic functions in $\C$ of the
form $u_j=\Rea P_j$ where $P_j\neq 0$ are polynomials. Suppose that
$u_1\neq u_2$.
Then the set $E=\{ z\in\C:u_1(z)=u_2(z)\}$ is a
proper real-algebraic
subset of $\bC$
without isolated points.
Apart from a finite set of ramification points, $E$ consists
of smooth curves. For every smooth point $z\in E$, we denote by
$J(z)$ the jump of the normal (to $E$)
derivative of the function $w=\max\{ u_1,u_2\}$
at the point $z$. This jump is always positive and the Riesz measure
$\mu_w$ is given by the formula
\begin{equation}\label{jump}
d\mu_w=\frac{J(z)}{2\pi}|dz|,
\end{equation}
which means that $\mu_w$ is supported by $E$ and has a density $J(z)/2\pi$
with respect to the length element $|dz|$ on $E$.

Now let $E_{i,j}=\{ z:u_i(z)=u_j(z)\geq u_k(z), 1\leq k\leq n\}$,
and $E=\cup E_{i,j}$
where the union is taken over all pairs $1\leq i,j\leq n$ for which
$u_i\neq u_j$. Then $E$ is a proper real semi-algebraic
subset of $\bC$,
and $\infty$ is not an isolated point of $E$. For the elementary properties
of semi-algebraic sets that we use here see, for example, \cite{BR,C}.
There exists $r_0>0$ such that $\Gamma=E\cap\{ r_0<|z|<\infty\}$
is a union of finitely many disjoint smooth simple curves,
$$\Gamma=\cup_{k=1}^m \Gamma_k.$$
This union coincides with the support of $\mu_{u^*}$ in $\{ z:r_0<|z|<\infty\}$.

Consider a point $z_0\in \Gamma$. Then $z_0\in\Gamma_k$ for some $k$.
As $\Gamma_k$ is a smooth curve, there is a neighborhood $D$ of $z_0$
which does not contain other curves $\Gamma_j,\; j\neq k$ and
which is divided by $\Gamma_k$ into two parts, $D_1$ and $D_2$.
Then there exist $i$ and $j$ such that $u^*(z)=u_i(z),\; z\in D_1$
and $u^*(z)=u_j(z),\; z\in D_2$, and $u^*(z)=\max\{ u_i(z),u_j(z)\},\;
z\in D$. So the restriction of the Riesz measure $\mu_{u^*}$ on
$D$ is supported by $\Gamma_k\cap D$ and has density
$J(z)/(2\pi)$
where 
$$|J(z)|=|\partial u_i/\partial n-\partial u_j/\partial n|(z)
=|\Del(u_i-u_j)|(z),$$
and $\partial/\partial n$ is the derivation in the direction
of a normal to $\Gamma_k$. Taking into account that $u_j=\Rea P_j$
where $P_j$ are polynomials, we conclude that there exist
positive numbers $c_k$ and $b_k$
such that
\begin{equation}
\label{jump1} J(z)/(2\pi)=(c_k+o(1))|z|^{b_k},\quad z\to\infty,\quad z\in\Gamma_k.\end{equation}
Let $b=\max_kb_k$, and among those curves $\Gamma_k$ for which
$b_k=b$ choose one with maximal $c_k$ (which we denote by $c_0$).
We denote this chosen curve by $\Gamma_0$ and fix it for the rest of the
proof. 
\vspace{.1in}

\noindent
{\bf Proposition 3.} {\em We have
$$b\leq\sigma\quad\mbox{and}\quad c_0\leq 3\cdot 4^{\sigma}K(n+1). $$
}

{\em Proof.}
We consider two cases.

Case 1. There is a sequence $z_n\to\infty,\; z_n\in\Gamma_0$ such that
$u_0(z_n)\leq u^*(z_n)$. Then (\ref{cond}) and Proposition 1 imply that
$$J(z_n)\leq (n+1)K|z_n|^\sigma,$$
and comparison with (\ref{jump1})
shows that $b\leq\sigma$ and $c_0\leq K(n+1)/(2\pi)$.

Case 2. $u_0(z)>u^*(z)$ for all sufficiently large $z\in\Gamma_0$.
Let $a$ be a point on $\Gamma_0$, $|a|>3r_0$, and $u_0(a)>u^*(a)$.
Let $B(a,R)$ be the largest
open disc centered at $a$ in which the inequality $u_0(z)>u^*(z)$
holds. Then
\begin{equation}
\label{R}
R\leq |a|+O(1),\quad a\to\infty
\end{equation}
because we assume that $f_0$ has zeros,
so $u_0(z_0)=-\infty$ for some $z_0$.

In $B(a,R)$ we consider the positive superharmonic function
$v=u_0-u^*$. Let us check that it satisfies the conditions of Lemma 2.
The existence of a point $z_1\in\partial B(a,R)$ with $v(z_1)=0$
follows from the definition of $B(a,R)$. The Riesz measure of $\mu_v$
is estimated using (\ref{jump}), (\ref{jump1}):
$$|\mu_v(B(a,R/2))|\geq |\mu_v(\Gamma_0\cap B(a,R/2))|\geq c_0R(|a|-R/2)^b.$$
Now Lemma 2 applied to $v$ in $B(a,R)$ implies that
\begin{equation}
\label{del}
|\Del v(z_1)|\geq (c_0/3)(|a|-R/2)^{b}.
\end{equation}
On the other hand (\ref{cond}) and Proposition 1 imply that
$$|\Del v(z_1)|\leq K(n+1)(|a|+R)^\sigma$$
Combining these two inequalities and taking (\ref{R}) into account,
we obtain $b\leq\sigma$ and $c_0\leq 3\cdot 4^{\sigma}K(n+1),$
as required.
\vspace{.1in}

We denote $$T^*(r)=\frac{1}{2\pi}\int_{-\pi}^\pi
u^*(re^{i\theta})d\theta -u^*(0);$$
This is the characteristic of the ``reduced curve''
$(f_1,\ldots,f_n)$. 
\vspace{.1in}

\noindent
{\bf Proposition 4.} 
$$T^*(r)\leq 6\cdot 4^\sigma K\frac{n(n+1)^2}{\sigma+1}.$$

{\em Proof.}
By Jensen's formula,
$$T^*(r)=\int_0^r\nu(t)\frac{dt}{t},$$
where $\nu(t)=\mu_{u^*}(\{ z:|z|\leq t\})$. The number of curves
$\Gamma_k$ supporting the Riesz measure of $u^*$ is
easily seen to
be at most $2n(n-1)(\sigma+1)$. The density of the Riesz measure
$\mu_{u^*}$ on each curve $\Gamma_k$ is given by (\ref{jump1}), where $c_k\leq c_0$
and $b_k\leq b$, and the parameters $c_0$ and $b$ are estimated in
Proposition 3. Combining all these data we obtain the result.
\vspace{.1in}

It remains to combine Propositions 2 and 4 
to obtain the final result.

{\em Purdue University, West Lafayette IN 47907 USA

}
\end{document}